\documentclass{ifacconf}

\usepackage{graphicx}
\usepackage{natbib}

\usepackage{amssymb,amsmath,amsfonts,cases}
\usepackage{soul}
\usepackage[dvipsnames]{xcolor}
\usepackage{algorithm}
\usepackage{algpseudocode}
\usepackage{algorithmicx}
\usepackage{enumitem}
\usepackage{subfig} % required for bibliography
\newcommand{\R}{\mathbb{R}}

\newcommand{\N}{\mathbb{N}}

\newcommand{\T}{^\top}
\newcommand{\inv}{^{-1}}

\newcommand{\argmin}{\mathop{\rm argmin}}

\newcommand{\subj}{\textnormal{subj.~to}}

\DeclareMathOperator{\col}{col}

\newtheorem{theorem}{Theorem}

\newtheorem{lemma}{Lemma}
\newtheorem{assumption}{Assumption}
\newtheorem{remark}{Remark}

\newcommand\oprocendsymbol{\hbox{$\blacksquare$}}
\newcommand\oprocend{\relax\ifmmode\else\unskip\hfill\fi\oprocendsymbol}

\newcommand{\cE}{\mathcal E}

\newcommand{\cG}{\mathcal G}
\newcommand{\cH}{\mathcal H}

\newcommand{\cL}{\mathcal L}

\newcommand{\cN}{\mathcal N}

\newcommand{\cV}{\mathcal V}

\renewcommand{\u}{\boldsymbol{u}}

\renewcommand{\int}{\mathbb{Z}}

\renewcommand{\int}{\mathbb{I}}

\newcommand{\norm}[1]{{\left\vert\kern-0.3ex\left\vert #1\right\vert\kern-0.3ex\right\vert}}
\newcommand{\mnorm}[1]{{\left\vert\kern-0.3ex\left\vert\kern-0.3ex\left\vert#1\right\vert\kern-0.3ex\right\vert\kern-0.3ex\right\vert}}

\newcommand{\g}{\boldsymbol{g}}
\newcommand{\s}{\boldsymbol{s}}
\newcommand{\x}{\boldsymbol{x}}

\newcommand{\z}{\boldsymbol{z}}

\renewcommand{\u}{\boldsymbol{u}}

\newcommand{\p}{\boldsymbol{p}}

\newcommand{\bzero}{\boldsymbol{0}}
\renewcommand{\int}{\mathbb{Z}}
\newcommand{\G}{\mathcal{G}}

\graphicspath{{figs/}}

\renewcommand{\norm}[1]{\left \|#1 \right \|}

\newcommand{\ud}{^}
\newcommand{\du}{_}

\newcommand{\iter}{t}
\newcommand{\iterp}{{\iter+1}}

\newcommand{\n}{n}
\newcommand{\m}{m}

\newcommand{\xt}{\x\ud\iter}
\newcommand{\xtp}{\x\ud{\iterp}}

\newcommand{\zt}{\z\ud\iter}
\newcommand{\ztp}{\z\ud{\iterp}}

\newcommand{\bz}{\bar{\z}}
\newcommand{\bzt}{\bz\ud\iter}

\newcommand{\pz}{\z_\perp} 
\newcommand{\pzt}{\pz\ud\iter} 
\newcommand{\pztp}{\pz\ud\iterp}

\newcommand{\lipp}{c}

\newcommand{\xit}{\x\ud\iter_{i}}
\newcommand{\xjt}{\x\ud\iter_{j}}
\newcommand{\xitp}{\x\ud{\iterp}_{i}}

\newcommand{\xstar}{x\ud\star}

\newcommand{\f}{f}

\newcommand{\J}{J}

\newcommand{\1}{\mathbf{1}}

\newcommand{\lt}{\lambda\ud\iter}
\newcommand{\ltp}{\lambda\ud\iterp}

\newcommand{\lit}{\lambda\ud\iter_i}
\newcommand{\ljt}{\lambda\ud\iter_j}
\newcommand{\litp}{\lambda\ud\iterp_i}

\newcommand{\zit}{\zt_i}
\newcommand{\zjt}{\zt_j}

\newcommand{\bx}{\bar{\x}}

\newcommand{\step}{\gamma}
\newcommand{\degi}{\text{deg}_i}

\newcommand{\lstar}{\lambda\ud\star}

\newcommand{\agg}{\alpha}

\newcommand{\al}{g}

\newcommand{\nb}{\bar{\nz}}

\newcommand{\pzeq}{\pz\ud{\text{eq}}}

\newcommand{\tz}{\tilde{\z}_\perp}

\newcommand{\Vz}{U}
\newcommand{\Vc}{W}

\newcommand{\nz}{p}

\newcommand{\tagg}{\hat{\agg}}
\newcommand{\taggi}{\tagg_i}
\newcommand{\taggj}{\tagg_j}

\newcommand{\ptagg}{\tagg_\perp}

\newcommand{\zij}{\z_{ij}}
\newcommand{\zji}{\z_{ji}}

\newcommand{\zijt}{\zij\ud\iter} 
\newcommand{\zjit}{\zji\ud\iter} 

\newcommand{\zijtp}{\zij\ud\iterp}

\newcommand{\np}{\nz_\perp}

\newcommand{\Tp}{T_\perp}

\newcommand{\reg}{\kappa}

\newcommand{\str}{\mu}

\newcommand{\dd}{\mathbf{d}}

\newcommand{\gp}{\g_\perp}

\newcommand{\chit}{\chi\ud\iter}
\newcommand{\chitp}{\chi\ud\iterp}

\newcommand{\chistar}{\chi\ud{\star}}

\newcommand{\aggx}{\agg^{\x}}
\newcommand{\aggl}{\agg^{\lambda}} 

\newcommand{\taggix}{\taggi^{\x}}
\newcommand{\taggil}{\taggi^{\lambda}} 

\newcommand{\taggjx}{\taggj^{\x}}
\newcommand{\taggjl}{\taggj^{\lambda}} 

\newcommand{\taggx}{\tagg^{\x}}
\newcommand{\taggl}{\tagg^{\lambda}} 

\newcommand{\ptaggx}{\ptagg^{\x}}
\newcommand{\ptaggl}{\ptagg^{\lambda}} 

\newcommand{\msg}{m}
\newcommand{\msgijt}{\msg_{ij}^\iter}
\newcommand{\msgjit}{\msg_{ji}^\iter}

\def\er/{Erd\H{o}s-R\'enyi}

\begin{document}

\begin{frontmatter}

	\title{Distributed Constraint-Coupled Optimization: Harnessing ADMM-consensus for robustness}
	%\thanks[footnoteinfo]{Work supported and funded by the Italian Ministry for Research
%and Education (MUR) under the National Recovery and Resilience Plan
%(NRRP)-MISSION 4 COMPONENT 2, Investment 1.3.}

	\author[First]{A.Mohamed Messilem} 
	\author[Second]{Guido Carnevale} 
	\author[Third]{Ruggero Carli}

	\address[First]{Department of Information Engineering (DEI), University of Padova, Padova, Italy
	(e-mail:
    mohamedabdelm.messilem@phd.unipd.it).}
	\address[Second]{Department of Electrical,  Electronic and Information Engineering,  Alma Mater Studiorum - Universita` di Bologna,  Bologna, Italy (e-mail: guido.carnevale@unibo.it)}
	\address[Third]{Department of Information Engineering (DEI), University of Padova, Padova, Italy
	(e-mail:carlirug@dei.unipd.it).}

	\begin{abstract}                % Abstract of 50--100 words
In this paper, we consider a network of agents that jointly aim to minimise the sum of local functions subject to coupling constraints involving all local variables.  To solve this problem, we propose a novel solution based on a primal-dual architecture. The algorithm is derived starting from an alternative definition of the Lagrangian function, and its convergence to the optimal solution is proved using recent advanced results in the theory of time-scale separation in nonlinear systems. The rate of convergence is shown to be linear under standard assumptions on the local cost functions.
Interestingly, the algorithm is amenable to a direct implementation to deal with asynchronous communication scenarios that may be corrupted by other non-idealities such as packet loss.
We numerically test the validity of our approach on a real-world application related to the provision of ancillary services in three-phase low-voltage microgrids.
	\end{abstract}


\end{frontmatter}

\section{Introduction}
The design of distributed optimization algorithms has received considerable attention in the last decade, driven by their wide range of applications in multi-agent systems, sensor networks, machine learning, and smart grids \cite{nedic2009distributed}. These algorithms allow agents networks to collaboratively solve optimization problems without relying on a central authority, making them particularly attractive for scalable, resilient, and privacy-preserving systems \cite{huo2021privacy}. 

In this paper, we are interested in problems where a network of agents aims to minimize the sum of local functions subject to coupling constraints. Such problems arise in many real-world applications, including the LASSO problems \cite{hastie2005elements} in machine learning, the power flow and
load control problems in smart grid \cite{alizadeh2012demand}, the network flow problem \cite{bertsekas1998network} and the coordinated transmission
design problem in communication networks \cite{shen2012distributed}, to name a few.
The presence of coupling constraints involving all local variables greatly complicates the design of algorithms, especially in a distributed setting where communication and coordination between agents is inherently challenging \cite{zhu2011distributed}.

Various works have focused on addressing the setup described above, with most solutions taking advantage of alternate direction method of multipliers (ADMM)-based approaches \cite{boyd2011distributed}. \cite{messilem2024distributed}. These approaches often combine gradient tracking with linear consensus steps and have demonstrated remarkable performance in synchronous environments under ideal conditions \cite{nedic2020distributed}. However, it is recognised that
making these strategies robust to asynchronous implementations and other non-idealities, such as packet loss or communication delays, requires consistent modifications of the algorithms, including additional variables and increasing the computational load.

In this paper, we propose a novel distributed primal-dual algorithm designed to address the challenges posed by coupling constraints. The derivation of the algorithm is obtained starting from an alternative definition of the Lagrangian function, and its convergence to the optimal solution is proved in the synchronous case, exploiting recent advances in time-scale separation theory. It is also shown that the convergence is linear under mild assumptions on the local cost functions. Interestingly, we illustrate how the algorithm can be easily adapted to work with asynchronous updates and is inherently robust to packet loss.
We test numerically our approach in real-world application related to the provision of Ancillary Services in Three-Phase Low-Voltage Microgrids.
%The paper is organized as follows. Section 2 formulates the problem, followed by the design of a distributed algorithm in Section 3. Section 4 provides an analysis of the proposed distributed algorithm, while Section 5 extends the framework to a robust distributed algorithm. Finally, Section 6 validates the effectiveness of our approach through numerical simulations.
\paragraph*{Notation}%
$\col(x_1,\dots,x_N)$ denotes the column stacking of the vectors $x_1,\dots,x_N$.
%
%Given $v := \col(v_1,\dots,v_n) \in \R^n$, 
%
In $\R^{m\times m}$, $I_m$ and $0_m$ are the identity and zero matrices. 
%
% The symbol $1_N$ denotes the vector of $N$ ones while $\1_{N,d} := 1_N \otimes I_d$, where $\otimes$ denotes the Kronecker product. 
The symbol $1_N$ denotes the vector of $N$.
Dimensions are omitted when they are clear from the context.
The symbol $\otimes$ denotes the Kronecker product. 
%
%Given $v \!\in\! \R^n$, we denote as $[v]_i$ its $i$-th component. %and as $\diag(v)$ the $n$-dimensional diagonal matrix whose $i$-th diagonal entry is $[v]_i$.
%
%Given $r > 0$, the symbol $\cB_r$ denotes the ball with radius $r$.
%
%Given a vector $x \in \R^{n}$ and a set $\Xv \subseteq \R^{n}$, $\Px{x}$ denotes the projection of $\x$ on $\Xv$, i.e., $\Px{x} := \argmin_{y \in \Xv}\norm{x - y}$.
%
%Given a dynamical system $\xtp = f(\xt)$ with $\xt \in \R^{\n}$ and $f: \R^\n \to \R^\n$, we say that the subset $\X \subseteq \R^\n$ is forward-invariant for the considered system if $\xt \in \X$ implies $\x\ud{\bar{\iter}} \in \X$ for all $\bar{\iter} \ge \iter$.
%
% Given $x \in \R^{n}$ and $\cS \subset \R^{\n}$, we define $\norm{x}_{\cS}:= \inf_{y \in \cS}\norm{x - y}$.
% %
% Given $f: \R^{n_1} \! \times \! \R^{n_2} \! \to \! \R$, we define $\nabla_1 f(x,y) \!:=\! \tfrac{\partial}{\partial s}f(s,y)|_{s = x}\T$ and $\nabla_2 f(x,y) \!:= \!\tfrac{\partial}{\partial s}f(x,s)|_{s = y}\T$. 
% %
% $\R_{+}^n$ is the positive orthant in $\R^n$.
Given a matrix $M \in \R^{n \times m}$, $\ker(M)$ denotes its kernel.

\section{Problem formulation}

Consider a network of $N$ agents modeled by an undirected connected graph $\mathcal{G}(\mathcal{V}, \mathcal{E})$, where $\mathcal{V}= \left\{1,\ldots, N\right\}$ is the set of nodes and $\mathcal{E}$ is the set of edges describing the admissible communications among the agents. 
The fact that $(i,j) \in \mathcal{E}$ means that $i$ and $j$ can communicate directly with each other. 
We denote by $\cN_i := \{j \in \cV \mid (i,j) \in \cE\}$ the set of neighbors of agent $i$ and with $\degi = |\cN_i|$ its degree.
Finally, we define $\dd := \sum_{i=1}^N\degi$, namely, $\dd := 2|\cE|$.

The $N$ agents cooperate to solve the optimization problem 
\begin{align} 
	\begin{split}\label{eq:constr_coupled}
		\min_{(x_1 \dots, x_N) \in \R^{\n}} \: & \: \sum_{i=1}^N \f_i(x_i)
		\\
		\subj \: & \: \sum_{i=1}^N A_ix_i = b,
		% % & \: \sum_{i=1}^N h_i(x_i) = 0\\
		% \: & \: \xag[i]\in \X[i], \forall i \in \set,
	\end{split}
\end{align}
where $f_i: \R^{\n_i} \to \R$, $i \in \mathcal{V}$, is the local cost function of agent
$i$ depending on the associated decision variable $x_i \in \R^{\n_i}$, and where the elements $A_i \in \R^{\m \times \n_i}$, $i\in \mathcal{V},$ and $b \in \R^{\m}$ model the coupled equality constraint.
The whole decision variable is denoted as
$x := \col(x_1,\dots,x_N) \in \R^{\n}$ with $\n := \sum_{i=1}^N
\n_i$.
For the sake of a compact notation, we introduce $\f: \R^{\n} \to \R$ and $A \in \R^{\m \times \n}$ defined as 
\begin{align*}
	\f(x) := \sum_{i=1}^N \f\du{i}(x\du{i}), \quad A := \begin{bmatrix}
		A\du{1}& \hdots& A\du{N}
	\end{bmatrix}.
\end{align*}
The following assumptions characterize $\f$ and $A$.
\begin{assumption}\label{ass:cost}
	There exists $\str > 0$ such that $\sum_{i=1}^N\f_i$ is $\str$-strongly convex.
	Moreover, there exists $\lipp > 0$ such that $\nabla\f_i$ is $\lipp$-Lipschitz continuous for all $i \in \cV$.\oprocend
\end{assumption}
\begin{assumption}\label{ass:full_row}
	The matrix $A$ is full-row rank.\oprocend
\end{assumption}
%
% We assume that the function $f_i$ is convex.

Let us now introduce the standard Lagrangian function $L: \R^{\n} \times \R^{\m} \to \R$ associated to problem~\eqref{eq:constr_coupled}, namely
\begin{align*}
	L(x, \bar{\lambda})= \sum_{i=1}^N \f_i(x_i) + \bar{\lambda}^T\left(\sum_{i=1}^N A_ix_i - b\right),
\end{align*}
where $\bar{\lambda} \in \R^m$ is the vector of Lagrange multipliers. 
%
% The dual problem of \eqref{eq:constr_coupled} is 
% \begin{equation}\label{eq:dual_problem}
% \max_{\bar{\lambda} \in \R^{\m}} \min_{x\in \R^{\n}} L(x, \bar{\lambda})     
% \end{equation}
%
We remark that Assumptions~\ref{ass:cost} and~\ref{ass:full_row} ensure that $L$ has a unique saddle-point $(\xstar,\lstar) \in \R^{\n} \times \R^{m}$, in which $\xstar$ is the unique solution to problem~\eqref{eq:constr_coupled} (see, e.g.~\cite{qu2018exponential}). 
% the existence of a unique solution $\xstar \in \R^{\n}$ to problem~\eqref{eq:constr_coupled} and an associated optimal multiplier $\lstar \in \R^{\m}$ being a solution to~\eqref{eq:dual_problem}.
% The following assumption guarantees that Problem in \eqref{eq:constr_coupled} and Problem in \eqref{eq:dual_problem} are well posed.
% \begin{assumption}
% Problem in \eqref{eq:constr_coupled} admits an optimal solution $x^*=\left[{x_1^*}^T \cdots {x_N^*}^T\right]^T$ and problem in \eqref{eq:dual_problem} admits an optimal solution $\lambda^*$.
% \end{assumption}
In a centralized setting, the above problems are addressed by resorting to dual-ascent algorithm or to an ADMM-like parallel algorithm (see \cite{bertsekas2015parallel}).
In view of distributed implementations, we introduce $N$ copies of the Lagrange multiplier $\bar{\lambda}$, that is, $\lambda_1, \ldots, \lambda_N \in \R^{\m}$, and we consider an alternative Lagrangian function $\cL: \R^{\n} \times \R^{N\m} \to \R$ defined as  
\begin{align}\label{eq:lagrangian}
    \cL(x,{\lambda}) &= \sum_{i=1}^N \f_i(x_i) + \frac{1}{N}\sum_{i=1}^N \lambda_i\T (Ax - b) 
	\notag\\
	&\hspace{.4cm}
	- \frac{\reg}{2} \lambda\T \left(I - \dfrac{1}{N}\1\1\T \right)\lambda,
\end{align}
where $\reg > 0$, $\lambda := \col(\lambda_1,\dots,\lambda_N) \, \in \R^{N\m}$, and $\1 := 1_N \otimes I_{\m}$.
%  and where we have compactly written the constraint by introducing $A \in \R^{\m \times \n}$ defined as 
% %
% \begin{align}
%     A := \begin{bmatrix}
%         A_1& \dots& A_N
%     \end{bmatrix}.
% \end{align}
% \GC{Potremmo anche dire che cercare punti di sella di $L$ (e poi di $\cL$) è equivalente a cercare punti di minimo di problema~\eqref{eq:constr_coupled}.}

Under Assumptions~\ref{ass:cost} and~\ref{ass:full_row}, for all $\reg > 0$, by observing the gradient of $\cL$ it is easy to see that the unique saddle-point of $\cL$ is $(\xstar,\1\lstar)$.
We then address problem~\eqref{eq:constr_coupled} by designing a distributed primal-dual method tailored to find the saddle-point of $\cL$.
More in detail, our aim is to get a distributed strategy that effectively works over networks that present the following challenges: (i) asynchronous operations of the agents and (ii) loss of peer-to-peer communications.
\begin{remark}	
In recent years, algorithms for solving problems in \eqref{eq:constr_coupled} have already been proposed in the literature.  
Most of them are ADMM-based strategies that include gradient tracking steps combined with linear consensus steps, (see \cite{falsone2020tracking,carli2019distributed}). These algorithms show remarkable performance when running in a synchronous scenario and with proper initializations of some variables. However, it is known that making these algorithms robust to asynchronous implementations and to other non-idealities, such as packet losses, requires consistently modifying the algorithms including additional variables and increasing the computational load.

\end{remark}

\section{Distributed Algorithm Design}

We now design a distributed primal-dual version to address problem~\eqref{eq:constr_coupled}.
By considering the Lagrangian function $\cL$ in~\eqref{eq:lagrangian}, agent $i$ would run 
\begin{subequations}\label{eq:centralized_primal_dual}
	\begin{align}
		\xitp &= \xit - \step \frac{\partial \cL(x,\lambda)}{\partial x_i} 
		\\
		\litp &= \lit + \step \frac{\partial \cL(x,\lambda)}{\partial \lambda_i}
 	\end{align}
\end{subequations}
that is,
\begin{subequations}\label{eq:centralized_primal_dual}
	\begin{align}
		\xitp &= \xit - \step\bigg(\nabla \f_i(\xit) + A_i\T \frac{1}{N}\sum_{j=1}^N\ljt\bigg)  
		\\
		\litp &= \lit + \step\bigg(\reg\bigg(\frac{1}{N}\sum_{j=1}^N\ljt - \lit\bigg) + \frac{1}{N}(A\xt - b)\bigg).
 	\end{align}
\end{subequations}
We note that algorithm~\eqref{eq:centralized_primal_dual} would violate our desired distributed paradigm since its execution would require the local knowledge of the global quantities $\frac{1}{N}\sum_{i=1}^N \lit$ and $\frac{1}{N}(A\xt - b)$.
For this reason, we now design a distributed version of~\eqref{eq:centralized_primal_dual} in which a concurrent auxiliary scheme provides local proxies of both $\frac{1}{N}\sum_{i=1}^N \lit$ and $\frac{1}{N}(A\xt - b)$ and iteratively improves them by running a distributed consensus-based dynamics.
To formally describe this concurrent scheme, given $\xt := \col(\xt\du{1},\dots,\xt\du{N}) \in \R^{\n}$ and $\lt := \col(\lt\du{1},\dots,\lt\du{N}) \in \R^{N\m}$, we compactly describe the global variable needed by~\eqref{eq:centralized_primal_dual} through the function $\agg: \R^{\n} \times \R^{\N\m} \to \R^{2\m}$ defined as
\begin{align}\label{eq:agg}
	\agg(x,\lambda) := \begin{bmatrix}
		\aggx(x) 
		\\
		\aggl(\lambda)
	\end{bmatrix} := \frac{1}{N}\sum_{i=1}^N\begin{bmatrix}
		A_ix_i - b 
		\\
		\lambda_i
	\end{bmatrix}.
\end{align}
Once this unavailable global quantity has been identified, in each agent $i \in \cV$, we introduce an auxiliary variable $\zijt \in \R^{2\m}$ for each neighbor $j \in \cN_i$ and replace $\agg(\xt,\lt)$ through a local proxy $\taggi(\col(\xit,\lit),\zit) \in \R^{2\m}$  only depending on the local variables $\xit$, $\lit$, and $\zit := \col(\zijt)_{j \in \cN\du{i}}$.
In particular, the local proxies and the auxiliary variables are iteratively updated according to consensus-ADMM (see~\cite{bastianello2022admm}), namely
 \begin{subequations}\label{eq:open_form_admm}
	\begin{align}
		\begin{bmatrix}
		    \taggix(\xit,\zit) 
            \\
            \taggil(\lit,\zit)	
		\end{bmatrix}
		&= \argmin_{\s_i \in \R^{2\m}}
		\bigg\{\tfrac{1}{2}\norm{\s_i - \begin{bmatrix}
			A_i\xit - b
            \\
			\lit
		\end{bmatrix}}^2 
		\notag\\
		&\hspace{1.7cm}
		- \s_i\T \sum_{j \in \cN_i}\zijt 
		+ \tfrac{\rho \degi}{2}\norm{\s_i}^2\bigg\}\!
		\label{eq:wit_cit}
		\\
		\zijtp &= (1 - \beta)\zijt + \beta\left(-\zjit + 2\rho\begin{bmatrix}
			\taggjx(\xjt,\zjt) 
            \\
            \taggjl(\ljt,\zjt)	
		\end{bmatrix}\right)
		\notag\\
		&\hspace{1cm}\forall j \in \cN_i,
	\end{align}
\end{subequations}
in which $\rho > 0$ and $\beta \in (0,1)$ are tuning parameters.
Since the cost function in~\eqref{eq:wit_cit} is quadratic, we can compute its minimum in closed form and accordingly rewrite~\eqref{eq:open_form_admm} as
\begin{subequations}\label{eq:ADMM_update}
	\begin{align}
		\begin{bmatrix}
		    \taggix(\xit,\zit) 
            \\
            \taggil(\lit,\zit)	
		\end{bmatrix} &= \dfrac{1}{1+\rho \degi}
        \left(
        \begin{bmatrix}    
			A_i\xit - b_i
			\\
			\lit
		\end{bmatrix} + \sum_{j\in\cN_i} \zijt \right)
		\\
		\zijtp &= (1 - \beta)\zijt + \beta\left(-\zjit + 2\rho \begin{bmatrix}
			\taggjx(\xjt,\zjt) 
            \\
            \taggjl(\ljt,\zjt)	
		\end{bmatrix}\right)
		\notag\\
		&\hspace{1cm}\forall j \in \cN_i.
	\end{align}
\end{subequations}
By interconnecting the consensus mechanism~\eqref{eq:ADMM_update} with a modified version of the optimization method~\eqref{eq:centralized_primal_dual} using $\taggix(\xit,\zit)$ and $\taggil(\lit,\zit)$ in place of $\frac{1}{N}\sum_{j=1}^N\ljt$ and $\frac{1}{N}(A\xt - b)$, respectively, we get the overall distributed algorithm reported in Algorithm~\ref{algo:algo}.
% %
% \begin{subequations}\label{eq:distributed_primal_dual}
% 	\begin{align}
% 		% \begin{bmatrix}
% 		%     \wit 
%         %     \\
%         %     \cit
% 		% \end{bmatrix} &= \dfrac{1}{1+\rho \degi}
%         % \left(
%         % \begin{bmatrix}
% 		%     \lit
%         %     \\
%         %     A_i\xit - b_i
% 		% \end{bmatrix} + \sum_{j\in\cN_i} \zijt \right)
% 		% \\
% 		\xitp &= \xit - \step\left(\nabla \f_i(\xit) + A_i\T \taggil(\xit,\zit)\right)  
% 		\\
% 		\litp &= \lit + \step\left(\reg\left(\taggil(\lit,\zit) - \lit\right) + \taggix(\xit,\zit)\right)
% 		\\
% 		\zijtp &= (1 - \beta)\zijt + \beta\left(-\zjit + 2\rho \begin{bmatrix}
% 		    \taggjx(\xjt,\zjt) 
%             \\
%             \taggjl(\ljt,\zjt)	
% 		\end{bmatrix}\right)
% 		\notag\\
% 		&\hspace{1cm}\forall j \in \cN_i.
%  	\end{align}
% \end{subequations}
%
\begin{algorithm}[H]
	\begin{algorithmic}
		\State \textbf{Initialization}: $\x_i^0 \in \R^{\n}$, $\lambda\du{i}^0 \in \R^{\m}$, $\z_{ij}^0 \in \R^{2\m}$ $\forall j \in \cN_i$.
		\For{$\iter=0, 1, \dots$}
		\vspace{.1cm}
			\State $\begin{bmatrix}
			\taggix(\xit,\zit) 
            \\[.1cm]
            \taggil(\lit,\zit)	
			\end{bmatrix} = \frac{1}{1+\rho \degi}\left(\begin{bmatrix}A\du{i}\xit - b\\\lit\end{bmatrix} + \sum_{j\in\cN_i} \zijt\right)$
			\vspace{.1cm}
			\State $
			\begin{bmatrix}
				\xitp
				\\ 
				\litp
			\end{bmatrix}= \begin{bmatrix}
				\xit
				\\ 
				\lit
			\end{bmatrix} 
			- 
			\step
			\begin{bmatrix}
				\nabla \f_i(\xit) + A_i\T \taggil(\lit,\zit) 
				\\
				\reg(\taggil(\lit,\zit) - \lit) \!+\! \taggix(\xit,\zit)
			\end{bmatrix}$
			\For{$j \in \cN_i$}
				\State $\msgijt = -\zijt + 2\rho\begin{bmatrix}
					\taggix(\xit,\zit) 
					\\
					\taggil(\lit,\zit)	
				\end{bmatrix}$
				\State transmit $\msgijt$ to $j$ and receive $\msgjit$ from $j$
				\State $\zijtp = (1-\alpha)\zijt + \alpha\msgjit$
			\EndFor
		\EndFor
	\end{algorithmic}
	\caption{Distributed Algorithm (Agent $i$)}
	\label{algo:algo}
\end{algorithm}
We note that Algorithm~\ref{algo:algo} is purely distributed since each agent $i$ can run it by sending $\msgijt = -\zijt + 2\rho\begin{bmatrix}
	\taggjx(\xjt,\zjt) 
	\\
	\taggjl(\ljt,\zjt)	
\end{bmatrix}$ only to its neighbors $j \in \cN_i$ and by receiving $\msgjit$ only from them.
The next theorem ensures that Algorithm~\ref{algo:algo} solves problem~\eqref{eq:constr_coupled} with a linear rate of convergence.
\begin{theorem}\label{th:convegence}
	Consider Algorithm~\ref{algo:algo} and let Assumptions~\ref{ass:cost} and~\ref{ass:full_row} hold.
	Then, for all $\reg, \rho > 0$ and $\beta \in (0,1)$, there exist $\bar{\step} > 0$ such that, for all $\step \in (0,\bar{\step})$, $(\x\du{i}\ud0,\lambda\du{i}\ud0,\z\du{i}\ud{0}) \in \R^{\n\du{i}} \times \R^{\m} \times \R^{2\m\degi}$, $i \in \cV$, it holds 
	\begin{align}
		\norm{\begin{bmatrix}\xt - \xstar\\\lt - \lambda^\star\end{bmatrix}} \leq a_1\exp(-a_2\iter),
	\end{align}
	for all $\iter \in \N$ and some $a_1, a_2 >0$.
	\oprocend
\end{theorem}
The proof of Theorem~\ref{th:convegence} is provided in Section~\ref{sec:proof}.
Specifically, it uses the systematic procedure developed in~\cite{carnevale2024unifying} to analyze distributed algorithms through timescale separation arguments.

\section{Distributed Algorithm Analysis} 
\label{sec:analysis}

This section demonstrates that Algorithm~\ref{algo:algo} satisfies the set of sufficient conditions established by the systematic procedure in~\cite{carnevale2024unifying} to prove the effectiveness of distributed algorithms.
Specifically, the procedure in~\cite{carnevale2024unifying} identifies a set of properties related to the optimization-oriented and consensus-oriented components of Algorithm~\ref{algo:algo}, considered separately.

More in detail, after the aggregate algorithm reformulation provided in Section~\ref{sec:reformulation}, Section~\ref{sec:opt} shows that the centralized inspiring method~\eqref{eq:centralized_primal_dual} is able to solve problem~\eqref{eq:constr_coupled} with the properties formally stated in~\cite[Ass.~2 and~3]{carnevale2024unifying}.
Then, Section~\ref{sec:cons} establishes that the consensus-oriented component meets the conditions formally outlined in~\cite[Ass.~4,~5, and~6]{carnevale2024unifying}.
Finally, in Section~\ref{sec:proof}, we integrate these preparatory results to provide the proof of Theorem~\ref{th:convegence}.

\subsection{Algorithm Reformulation}
\label{sec:reformulation}

First of all, we provide the aggregate description arising from the local updates reported in Algorithm~\ref{algo:algo}.
In detail, given $\zt := \col(\zt\du{1},\dots,\zt\du{N}) \in \R^{2\m\dd}$, the aggregate evolution of Algorithm~\ref{algo:algo} reads as
\begin{subequations}\label{eq:distributed_primadual_aggregate}
	\begin{align}
		\xtp &= \xt - \step\left(\nabla\f(\xt) + A\T\taggl(\lt,\zt)\right)\label{eq:distributed_primadual_aggregate_x}
		\\
		\ltp &= \lt + \step(\reg(\taggl(\lt,\zt) - \lt) + \taggx(\lt,\zt))\label{eq:distributed_primadual_aggregate_l}
		\\
		\ztp &= T\zt + \g(\xt,\lt),\label{eq:distributed_primadual_aggregate_z}
	\end{align}
\end{subequations}
in which $T \in \R^{2\m\dd \times 2\m\dd}$, $\taggl: \R^{N\m} \times \R^{2\m\dd} \to \R^{N\m}$, $\taggx: \R^{N\m} \times \R^{2\m\dd} \to \R^{N\m}$, and $\g: \R^{\n} \times \R^{N\m} \to \R^{2\m\dd}$ are defined as 
\begin{align*}
	\taggl(\lambda,z)  
	&:= 
	\begin{bmatrix}
		\taggl\du{1}(\lambda\du{1},z\du{1})\T&\hdots&\taggl\du{N}(\lambda\du{1},z\du{1})\T
	\end{bmatrix}\T
    \\
    \taggx(x,z) 
	&:= 
	\begin{bmatrix}
		\taggl\du{1}(x\du{1},z\du{1})\T&\hdots&\taggl\du{N}(x\du{1},z\du{1})\T
	\end{bmatrix}\T
	\\
	T &:= I_{2\dd\m} - \beta(I_{2\dd\m} + S - 2\rho SD H D\T) 
	\\
	\g(x,\lambda) &:= 2\beta\rho SD H\col(A_1\xt_1 - b,\lt_1,\dots,A_N\xt_N - b,\lt_N),
\end{align*}
where $S \in \R^{2\m\dd \times 2\m\dd}$ is the permutation matrix that swaps the $ij$-th element with the $ji$-th one, while the matrices $D \in \R^{2\m\dd \times 2N\m}$ and $H \in \R^{2N\m \times 2N\m}$ read as
%

% 	\Ax &:=\begin{bmatrix}
% 			\one_{\degree_1,\n}
% 					\\
% 					0_{\degree_1,\n}
% 			\\
% 			&\diagentry{\xddots}
% 			\\
% 			&&
% 			\one_{\degree_N,\n}
% 			\\
% 			% &&&&
% 			%&
% 			&&0_{\degree_N,\n}
% 		\end{bmatrix} 
% %\\
% \!\!,\hspace{.08cm}
% 	\An := \begin{bmatrix}
% 			0_{\degree_1,\n}
% 			\\
% 			\one_{\degree_1,\n}
% 			\\
% 			&\diagentry{\xddots}
% 			\\
% 			&&0_{\degree_N,\n}
% 			\\		
% 			%&
% 			% \\
% 			% &&&&
% 			&&\one_{\degree_N,\n}
% 		\end{bmatrix}
% 		\\
	% \end{align*}
	% \begin{align*}
%
\begin{align*}
    D 
    &:=
    \begin{bmatrix}
    1_{\deg_1} \otimes I_{2m}
    \\
    &\ddots
    \\
    &&1_{\deg_N} \otimes I_{2m}
    \end{bmatrix}
    \\
    H 
    &:= 
    \begin{bmatrix}
     \frac{1}{1+\rho\deg_1}I_{2m}
     \\
    &\ddots
    \\
    &&
    \frac{1}{1+\rho\deg_N}I_{2m}
    \end{bmatrix}
\end{align*}
% \begin{align*}
%     D 
%     &:=
%     \begin{bmatrix}
%     1_{\deg_1} \otimes I_{2\m}
%     \\
%     &\diagentry{\xddots}
%     \\
%     &&1_{\deg_N} \otimes I_{2\m}
%     \end{bmatrix}
%     \\
%     H 
%     &:= 
%     \begin{bmatrix}
%      \frac{1}{1+\rho\deg_1}I_{2\m}
%      \\
%     &\diagentry{\xddots}
%     \\
%     &&
%     \frac{1}{1+\rho\deg_N}I_{2\m}
%     \end{bmatrix}
% \end{align*}
		% \\
		% \cH &:= \begin{bmatrix}
		% 	\frac{1}{1+\rho\degree_1}I_{2\n}
		% 			\\
		% 			&\diagentry{\xddots}
		% 			\\
		% 			&&
		% 			\frac{1}{1+\rho\degree_N}I_{2\n}
		% \end{bmatrix}.
		%\\
		%\Tar &:= (1 - \alpha)I - \alpha P + 2\alpha\rho PA\cH  A\T.

%
With the aggregate description~\eqref{eq:distributed_primadual_aggregate} at hand, we proceed by analyzing it by following the steps outlined in~\cite{carnevale2024unifying}.

\subsection{Optimization-Oriented Part Analysis}
\label{sec:opt}

In this section, we focus on the convergence properties of the inspiring, centralized algorithm~\eqref{eq:centralized_primal_dual} in solving problem~\eqref{eq:constr_coupled}.
%
% Indeed, in light of the perfect reconstruction condition~\eqref{eq:agg_reconstruction}, we note that subsystem~\eqref{eq:distributed_primadual_aggregate_tranformed_x}-\eqref{eq:distributed_primadual_aggregate_tranformed_tranformed_l} coincides with~\eqref{eq:centralized_primal_dual} in the case in which $\pzt = \pzeq(\xt,\lt)$ for all $\iter \in \N$.
%
To this end, we use the local update~\eqref{eq:centralized_primal_dual} to obtain its aggregate form
\begin{align}
	\chitp &= \chit + \step\al(\chit),\label{eq:compact_centralized}
\end{align}
where $\chit := \col(\xt,\lt) \in \R^{(\n+N\m)}$ and $\al: \R^{(\n+N\m)} \to \R^{(\n+N\m)}$ is defined as
\begin{align}
	\al(\chi) := 
	\begin{bmatrix}
		-\nabla \f(\x) - \frac{1}{N}A\T\1\T\lambda
		\\
		\frac{1}{N}\1A\x - b -\reg\J\lambda
	\end{bmatrix},\label{eq:al}
\end{align}
in which we used $\J := I_{N\m} - \frac{1}{N}\1\1\T$ and, with a slight abuse of notation, $\chi := \col(\x,\lambda)$ in the right-hand side.
The next lemma ensures that~\eqref{eq:compact_centralized} satisfies the conditions established by~\cite[Ass.~1, Ass.~2, and~3]{carnevale2024unifying}.
\begin{lemma}\label{lemma:centralized_method}
	The point $(\xstar,\1\lstar)$ is the unique equilibrium point of~\eqref{eq:compact_centralized}.
	Moreover, there exists a differentiable function $\Vc: \R^{(\n+N\m)} \to \R$ such that
	\begin{subequations}\label{eq:Vc_conditions}
		\begin{align}
			\underbar{c} \norm{\chi - \chistar}^2 \leq \Vc(\chi) &\leq \bar{c} \norm{\chi - \chistar}^2
			\label{eq:V_quadratic}   
			\\ 
			\nabla\Vc(\chi)\T \al(\chi) &\leq - c_1\norm{\chi - \chistar}^2
			\label{eq:V_decreasing}
			\\
			\norm{\nabla\Vc(\chi) - \nabla\Vc(\chi^\prime)} &\leq c_2\norm{\chi - \chi^\prime},
			\label{eq:Lipscitz_Vc} 
			\\
			\norm{\nabla \Vc(\chi)} &\leq c_3\norm{\al(\chi)},
			\label{eq:nabla_Vc_al}
 		\end{align}
	\end{subequations}
	for all $\chi \in \R^{(\n+N\m)}$ and some $\underbar{c},\bar{c},c_1,c_2,c_3 > 0$.\oprocend
\end{lemma}
The proof of lemma~\ref{lemma:centralized_method} is provided in Appendix~\ref{sec:proof_centralized}.
Roughly speaking, the conditions in~\eqref{eq:Vc_conditions} provide the ingredients to ensure that $\Vc$ is a Lyapunov function for system~\eqref{eq:compact_centralized} and the point $(\xstar,\lstar)$.
More in detail, the conditions in~\eqref{eq:Vc_conditions} ensure that, with sufficiently small $\step$, global exponential stability of $(\xstar,\lstar)$ for system~\eqref{eq:compact_centralized} can be achieved.

\subsection{Consensus-Oriented Part Analysis} 
\label{sec:cons}

Now, we need to show that the consensus-oriented part~\eqref{eq:distributed_primadual_aggregate_z} satisfies the conditions formalized by~\cite[Ass.~4,~5, and~6]{carnevale2024unifying}.

Specifically, in~\cite[Ass.~4 and~5]{carnevale2024unifying}, it is required that~\eqref{eq:distributed_primadual_aggregate_z} admits an orthogonal decomposition that allow us to isolate a part having equilibria parametrized in $(\xt,\lt)$ in which the global information $\agg(\xt,\lt)$ (cf.~\eqref{eq:agg}) is reconstructed through the proxies' vector $\tagg(\xt,\lt,\zt)$.

In order to find this decomposition, we follow the discussion provided in~\cite[Section~IV]{carnevale2025admm}.
In particular, we recall that since $\cG$ is connected, then $T$ has some simple eigenvalues equal to $1$ and all the remaining ones strictly inside the unit circle, see~\cite{bastianello2020asynchronous} for further details.
Let $\nb \in \N$ be the dimension of the subspace spanned by the unitary eigenvalues of $T$ and $B \in \R^{2\m\dd \times \nb}$ be the matrix whose columns are given by the eigenvectors associated to the unitary eigenvalues of $T$.
Analogously, given $\np = 2\m\dd - \nb$, let $M \in \R^{2\m\dd \times \np}$ be the matrix such that $M\T B = 0$ and $M\T M = I_{\np}$.
We now employ the matrices $M$ and $B$ to introduce the change of coordinates 
\begin{align}
	\z \longmapsto \begin{bmatrix}
		\bz 
		\\
		\pz
	\end{bmatrix} := \begin{bmatrix}
		B\T 
		\\
		M\T
	\end{bmatrix}\z.\label{eq:change_of_coordinates}
\end{align}
The discussion in~\cite[Section~IV]{carnevale2025admm} ensures that $\tagg$ does not depend on the new coordinate $\bz$, namely, there exists $\ptagg: \R^{\n} \times \R^{N\m} \times \R^{\np} \to \R^{2N\m}$ such that
\begin{align}\label{eq:orthogonality}
	\tagg(x,\lambda,\lambda,z) = \ptagg(x,\lambda,M\T\z) := 
	\begin{bmatrix}
		\ptaggx(x,M\T z)
		\\
		\ptaggl(\lambda,M\T z) 
	\end{bmatrix},
\end{align}
for all $x \in \R^{\n}$, $\lambda \in \R^{N\m}$, and $z \in \R^{2\m\dd}$, in which we further introduced $\ptaggx: \R^{\n} \times \R^{\np} \to \R^{N\m}$ and $\ptaggl: \R^{N\m} \times \R^{\np} \to \R^{N\m}$.
The orthogonality condition~\eqref{eq:orthogonality} allows us to disregard the evolution of $\bzt$ and focus only on $\pzt$.
By using~\eqref{eq:distributed_primadual_aggregate_z} and~\eqref{eq:change_of_coordinates}, $\pzt$ evolves according to
%
%\begin{subequations}\label{eq:distributed_primadual_aggregate_tranformed}
	\begin{align}
		% \xtp &= \xt - \step\left(\nabla\f(\xt) + A\T\ptaggl(\lt,\pzt) \right)\label{eq:distributed_primadual_aggregate_tranformed_x}
		% \\
		% \ltp &= \lt \! + \! \step(\reg(\ptaggl(\lt,M\pz) - \lt) \! + \! \ptaggx(\xt,\pzt))\label{eq:distributed_primadual_aggregate_tranformed_tranformed_l}
		% \\
		\pztp &= \Tp\pzt + \gp(\xt,\lt),\label{eq:distributed_primadual_aggregate_tranformed_pz}
	\end{align}
%\end{subequations}
%
in which $\Tp := M\T TM$ and $\gp(\xt,\lt) := M\T\g(\xt,\lt)$.
Note that system~\eqref{eq:distributed_primadual_aggregate_tranformed_pz} has equilibria parametrized in $(\xt,\lt)$ through the function $\pzeq: \R^{(\n+N\m)} \to \R^{\np}$ defined as 
\begin{align}
	\pzeq(x,\lambda) := (I_{\np}-\Tp)\inv\gp(x,\lambda).\label{eq:pzeq}
\end{align}
Indeed, we recall that $\Tp$ is the portion of the matrix $T$ obtained by removing its unitary eigenvalues of $T$. 
Thus, $\Tp$ is Schur and, in turn, the matrix $(I_{\np}-\Tp)$ in~\eqref{eq:pzeq} is invertible.
Moreover, in~\cite[Section~IV]{carnevale2025admm}, it is shown that 
%
%\begin{subequations}\label{eq:properties_pzeq}
	\begin{align}
		\ptagg(x,\lambda,\pzeq(x,\lambda)) &= \1\agg(x,\mathbf{\lambda}),
		\label{eq:agg_reconstruction}
		% \\
		% \pzeq(x) &= \Tp\pzeq(x) + \gp(x,\mathbf{\lambda}),
		% \label{eq:equilibrium_pzeq}
	\end{align}
%\end{subequations}
%
for all $(x,\mathbf{\lambda}) \in \R^{\n} \times \R^{N\m}$.
Namely, with $\pz = \pzeq(x,\lambda)$ perfect reconstruction of $\agg(x,\lambda)$ is achieved with the proxies' vector $\ptagg$.

In the ideal case in which $(\xt,\lt) = (x,\lambda)$ would be (arbitrarily) fixed, \cite[Ass.~6]{carnevale2024unifying} requires that $\pzeq(\xt,\lt)$ be globally exponentially stable for~\eqref{eq:distributed_primadual_aggregate_tranformed_pz}, uniformly in $(\xt,\lt) = (x,\lambda)$.
%
%
% We now focus on the stability properties of the equilibrium point $\pzeq(x,\lambda)$ for~\eqref{eq:distributed_primadual_aggregate_tranformed_pz}.
% %
% Since $\Tp$ has been obtained by removing the marginal stable part of $T$, we note the matrix $\Tp$ is Schur.
% %
% This property ensures global exponential stability of $\pzeq(x,\mathbf{\lambda})$ uniformly in $(x,\mathbf{\lambda})$ for~\eqref{eq:distributed_primadual_aggregate_tranformed_pz} in the ideal case in which $\xt$ and $\lt$ are fixed, namely, in the case in which $(\xt,\lt) = (x,\mathbf{\lambda})$ for all $\iter \in \N$.
%
We conclude this part by reporting the result~\cite[Lemma~IV.3]{carnevale2025admm} which formally provides a Lyapunov function ensuring the mentioned stability properties.
It is obtained by observing that~\eqref{eq:distributed_primadual_aggregate_tranformed_pz} is a linear system with Schur state matrix $\Tp$.
\begin{lemma}\label{lemma:lyapunov_pz}[\cite[Lemma~IV.3]{carnevale2025admm}]
	There exists $U: \R^{\np} \to \R$ such that 
	\begin{subequations}\label{eq:Vz_theorem_distributed}
		\begin{align}
			b_1\norm{\tz}^2 \leq \Vz(\tz) &\leq b_2\norm{\tz}^2\label{eq:Vz_quadratic_distributed}
			\\
			\Vz(\Tp\tz) -  \Vz(\tz) &\leq - b_3\norm{\tz}^2 
			\\
			|\Vz(\tz) - \Vz(\tz^\prime)| &\leq b_4\norm{\tz - \tz^\prime}(\norm{\tz} + \norm{\tz^\prime}),
		\end{align}
	\end{subequations}
	for all $(x,\lambda) \in \R^{\n} \times \R^{N\m}$, $\tz, \tz^\prime \in \R^{\np}$, and some $b_1, b_2, b_3, b_4 > 0$.\oprocend
\end{lemma}

\subsection{Proof of Theorem~\ref{th:convegence}}
\label{sec:proof}

We now conclude the proof of Theorem~\ref{th:convegence} by relying on~\cite[Th.~1]{carnevale2024unifying}.
To apply it, we need that the following conditions hold true.

First, by Lemma~\ref{lemma:centralized_method} and the Lipschitz continuity of the gradients $\nabla \f\du{i}$ (cf. Assumption~\ref{ass:cost}), we guarantee that the inspiring, centralized optimization method~\eqref{eq:centralized_primal_dual} is able to solve problem~\eqref{eq:constr_coupled} with the conditions formally prescribed by~\cite[Ass.~1,~2, and 3]{carnevale2024unifying}.

Second, the change of variables~\eqref{eq:change_of_coordinates}, the orthogonality condition~\eqref{eq:orthogonality}, the equilibrium function $\pzeq$ (cf.~\eqref{eq:pzeq}), reconstruction property~\eqref{eq:agg_reconstruction}, and the Lyapunov function $\Vz$ (cf. Lemma~\ref{lemma:lyapunov_pz}) ensure that the consensus mechanism~\eqref{eq:distributed_primadual_aggregate_z} satisfies the conditions formally established by~\cite[Ass.~3,~4, and~5]{carnevale2024unifying}.

Therefore, we are in position to apply~\cite[Th.~1]{carnevale2024unifying}.
Specifically, once we have chosen $\rho, \beta$, and $\reg > 0$, the result~\cite[Th.~1]{carnevale2024unifying} ensures there exists $\bar{\step} > 0$ such that, for all $\step \in (0,\bar{\step})$, the trajectories of the distributed algorithm~\eqref{eq:distributed_primadual_aggregate} satisfies
\begin{align*}
    \norm{\begin{bmatrix}
        \xt - \xstar 
        \\
        \lt - \1\lstar 
        \\
        \pzt - \pzeq(\xt,\lt)
    \end{bmatrix}} \leq a_0\norm{\begin{bmatrix}
        \x\ud0 - \xstar 
        \\
        \lambda\ud0 - \1\lstar 
        \\
        \p\ud0 - \pzeq(\x\ud0,\lambda\ud0)
    \end{bmatrix}}\exp(-a_2\iter),
\end{align*}
for all $\iter \in \N$, $(\x\ud0,\lambda\ud0,\pz\ud0)\in\R^{\n} \times \R^{N\m} \times \R^{\np}$, and some $a_0, a_2 > 0$, and the proof concludes.

\section{Robust Distributed Algorithm}
\label{sec:robust_algorithm} 

In this section, we provide a robust version of Algorithm~\ref{algo:algo} tailored for scenarios in which the network is subject to asynchronous agents and packet losses.
\begin{algorithm}[H]
	\begin{algorithmic}
		\State \textbf{Initialization}: $\x_i^0 \in \R^{\n}, \lambda\du{i}^0 \in \R^{\m}, \z_{ij}^0 \in \R^{2\m} \hspace{.1cm} \forall j \in \cN_i$.
		\For{$\iter=0, 1, \dots$}
		\If{active}
		\vspace{.1cm}
		\State $\begin{bmatrix}
			\taggjx(\xjt,\zjt) 
            \\[.1cm]
            \taggjl(\ljt,\zjt)	
			\end{bmatrix} = \frac{1}{1+\rho \degi}\left(\begin{bmatrix}A\du{i}\xit - b\\\lit\end{bmatrix} + \sum_{j\in\cN_i} \zijt\right)$
			\vspace{.1cm}
		\State 
		$\begin{bmatrix}
			\xitp
			\\ 
			\litp
		\end{bmatrix} \!=\! \begin{bmatrix}
			\xit
			\\ 
			\lit
		\end{bmatrix} 
		-
		\step
		\begin{bmatrix}
			\nabla \f_i(\xit) + A_i\T \taggil(\xit,\zit)  
			\\
			\reg(\taggil(\lit,\zit) - \lit) \!+\! \taggix(\xit,\zit)
		\end{bmatrix}$
			\For{$j \in \cN_i$}
						\State $\msgijt = -\zijt + 2\rho\begin{bmatrix}
							\taggjx(\xjt,\zjt) 
							\\[.1cm]
							\taggjl(\ljt,\zjt)	
						\end{bmatrix}$
						\State transmit $\msgijt$ to $j$ and receive $\msgjit$ to $j$
						\If{$\msgjit$ is received}
							\State $\zijtp = (1-\alpha)\zijt + \alpha\msgjit$
						\EndIf
			\EndFor
		\EndIf
		\EndFor
	\end{algorithmic}
	\caption{Robust Distributed Algorithm (Agent $i$)}
	\label{algo:ralgo}
\end{algorithm}
In Algorithm~\ref{algo:ralgo}, we note that each variable $\zijt$ is updated only if, at iteration $\iter$, the message from agent $j$ has been received by agent $i$.
Furthermore, as expected from the above discussion, the additional condition for performing such an update is that agent $i$ is active at iteration $\iter$.

\section{Numerical Simulations} 
\label{sec:numerical_simulations}
In this section we numerically validate our approach in a framework related to the provision of ancillary services in three-phase low-voltage microgrids.
The microgrid under consideration includes (i) controllable units that are distributed electronic power converters (EPCs) serving as interfaces between local energy resources, such as renewables or energy storage systems, and the grid, and (ii) uncontrollable
units, such as loads (e.g., households, offices). 
 In addition, the microgrid is connected
to the main grid at the point of common coupling (PCC).

We consider the control problem proposed in \cite{lauri2024ancillary} where the goal is to compute the optimal provision of reactive power, unbalanced currents, and zero sequence powers by distributed EPCs to minimize the converter power losses; as a consequence, this allows the quality of the power to be enhanced at the connection to the upstream grid.

%Numerical tests validate the proposed method for a network comprising four agents: three controllable EPCs and the grid, along with non-controllable units such as renewable sources and loads. The method formulates a control problem to provide ancillary services, such as reactive power compensation and unbalanced current management, while minimizing power losses and avoiding constraints on active power flows \cite{lauri2024ancillary}. The algorithm's operation is demonstrated using a cost function that includes converter power losses relative to their operating points, ensuring efficient and stable system performance.
Let $\x_i=\left[I^+_{d,i},I^+_{q,i},I^-_{d,i},I^-_{q,i},I^0_{d,i},I^0_{q,i}\right]$, be the vector with the dq components of the positive sequence, negative sequence, and zero sequence currents of the $i$-th controllable converter \footnote{In these problems $I^+_{d,i}$ is typically assumed to be constant and it is not considered as an optimization variable. Indeed, the goal is to minimize the power losses without changing the active power flows which are dictated by the needs of the grid and in the literature it is known that the active power exchange is approximately determined
by the direct component of the exchanged current.}. 
Moreover, let $\u_i\triangleq \left[\tilde{I}^+_{d,i},\tilde{I}^+_{q,i},\tilde{I}^-_{d,i},\tilde{I}^-_{q,i},\tilde{I}^0_{d,i},\tilde{I}^0_{q,i}\right]$, be the vector containing the quantities of the $i$-th non controllable unit ($\tilde{\cdot}$ refers to a non controllable quantity).
%Let ${\x \triangleq \left[\x_0,\ldots,\x_i,\ldots,\x_{N_{x}}\right]^\top}$ be the column vector of all optimization variables, \\with
%${\x_i=\left[I^+_{d,i},I^+_{q,i},I^-_{d,i},I^-_{q,i},I^0_{d,i},I^0_{q,i}\right]\in\mathbb{R}^{1\times6}}$, that is, the rms value of the reactive current of the positive sequence, of the negative sequence $d$ and $q$ currents, the zero sequence $d$ and $q$. 
%The vectors containing quantities referring to uncontrollable units are defined as ${\u_i\triangleq \left[\tilde{I}^+_{d,i},\tilde{I}^+_{q,i},\tilde{I}^-_{d,i},\tilde{I}^-_{q,i},\tilde{I}^0_{d,i},\tilde{I}^0_{q,i}\right]\in\mathbb{R}^{1\times 6}}$, where $\tilde{\cdot}$ denotes a noncontrollable quantity.

Assume to have $N_x$ controllable units. Then, the considered cost function accounts for converter losses and the power losses at the PCC interface as:
\begin{equation}
    F\left(\bx\right) = P^0_\mathit{loss}\left(\x_0\right) + \sum_{i=1}^{N_x}P^i_\mathit{loss}\left(\x_i\right)
    \label{eq:cost}
\end{equation}
where $P^0_\mathit{loss}\left(\x_0\right)$ is the term accounting for PCC interface losses ($i = 0$ refers to the quantities at the interface with
the main grid), while $P^i_\mathit{loss}\left(\x_i\right),~i=1,\ldots,N_x$ is the power loss of the $i$-th controllable converter.
Assuming a grid resistance value $R_g$, $P^0_{loss}(x_0)$
can be estimated by 
\begin{equation}
    P^0_\mathit{loss}\left(\x_0\right) = 3R_g\left\lVert\x_0\right\lVert^2
\end{equation}
while the losses on the converters can be modeled
as
 \begin{equation}
   P_\mathit{loss}^i\left(\x_i\right) = a_i\cdot\vert|\x_i\vert|^2 + b_i\cdot\vert|\x_i\vert| + c_i,
   \label{eq:ploss}
\end{equation}
where $a_i, b_i, c_i$ are suitable parameters.
%The microgrid under consideration includes uncontrollable units, such as loads (e.g., households, offices) and generators (e.g., wind turbines, photovoltaic panels), along with controllable units. In addition, the microgrid is connected to the main grid at the point of common coupling (PCC). 
According to Kirchhoff’s current laws, the combined power of the grid and controllable units must balance the net power demand of the uncontrollable units denoted as $\g$, i.e.,
\begin{equation}
    \g^\top = \sum_{i=1}^{N_u}\u_i = \x_0 + \sum_{i=1}^{N_x}\x_i
    \label{eq:constraint}
\end{equation}
where $N_u$ is the number of non controllable units.
Thus, our optimization problem can be written in a fully decoupled way as follows:
\begin{equation}
    \begin{aligned}
        \min_{\x_0, \ldots, x_{N_x}}& \sum_{i=0}^{N_x}f_i(\x_i)\\
        \text{s.t.} &\sum_{i=0}^{N_x}   \x_i  - \g = \bzero;  
    \end{aligned}
    \label{eq:problem_dis}
\end{equation}

The simulation was carried out considering three different cases (centralized, synchronous distributed, and asynchronous distributed). The stopping criterion for distributed algorithms is based on the distance from the optimal solution $d^k = \sum_{i=0}^{N_x}||\x_i^k - \x_i^*||_2^2$, where $x^k$ is the solution in iteration $k$ and $x^*$ is the optimal solution.
\begin{figure}[h]
    \centering
    \includegraphics[width=0.45\textwidth]{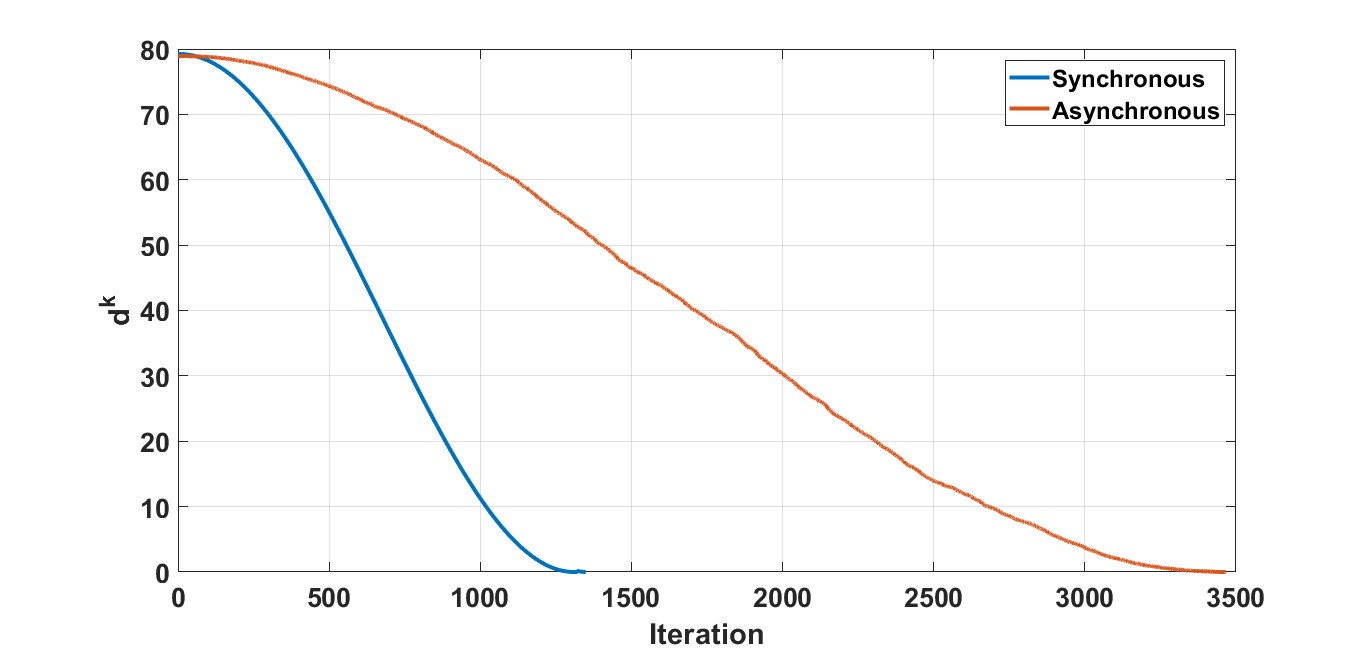}
    \caption{Numerical simulation results, comparing Synchronous and asynchronous approach.}
    \label{fig1}
\end{figure}

Figure \ref{fig1} illustrates the number of iterations required for synchronous and asynchronous distributed approaches to converge to the optimal solution. The plot demonstrates the algorithm's ability to converge to the optimal solution. Both methods achieve the same result. Indeed, the asynchronous approach requires a higher number of iterations to converge compared to the synchronous method, which aligns with expectations.

\section{Conclusions}

\label{sec:conclusions}
In this paper, we considered a distributed optimization problem to be solved over a network of agents with coupling constraints. Each agent represents a unit in the network and must cooperate with others to solve the problem without a centralized coordinator. We proposed a novel asynchronous distributed primal-dual algorithm to address the challenges posed by asynchrony and communication failures, providing theoretical guarantees on its convergence. Numerical simulations demonstrate the effectiveness of the proposed approach, showcasing its robustness and scalability in asynchronous scenarios. 

\appendix 

\section{Proof of Lemma~\ref{lemma:centralized_method}}
\label{sec:proof_centralized}

As for the first claim, it trivially follows by the fact that $\chistar := \col(\xstar,\1\lstar)$ is the unique saddle-point of the Lagrangian function $\cL$ (cf.~\eqref{eq:lagrangian}).

The remainder of the proof is inspired by the one in~\cite[Lemma~2]{qu2018exponential} in which, however, the considered Lagrangian function does not contain the term $-\frac{\reg}{2}\lambda\T\left(I - \frac{1}{N}\1\1\T\right)\lambda$.

Moreover, in~\cite[Lemma~2]{qu2018exponential}, a continuous-time scenario is considered.

We note that since $\chistar$ is a saddle-point of the Lagrangian function $\cL$ (cf.~\eqref{eq:lagrangian}), we can compactly rewrite the vector field $\al$ (cf.~\eqref{eq:al}) as

\begin{align}
	\al(\chi) := 
	\G(\chi,\chistar)(\chi - \chistar),\label{eq:G_al}
\end{align}

where $\G\!:\!\R^{(\n\!+\!N\m)} \!\times\! \R^{(\n+N\m)} \!\to\! \R^{(\n+N\m) \times (\n+N\m)}$ reads as

\begin{align}
	\G(\chi,\chistar) := 
	\begin{bmatrix}
		-\cH(\x,\xstar)& -\frac{1}{N}A\T\1\T
		\\
		\frac{1}{N}\1A& -\reg\J
	\end{bmatrix},\label{eq:G}
\end{align}

in which $\chi :=\col(\x,\lambda)$ and we also used the Mean Value Theorem to write $\nabla f(\x) - \nabla f(\xstar) = \cH(\x,\xstar)(\x - \xstar)$ where $\cH(\x,\xstar) := \nabla^2 \f(\bar{x}(\x,\xstar))$ for some $\bar{x}(\x,\xstar) \in \R^{\n}$.

Now, we pick $W: \R^{(\n+N\m)} \to \R$ defined as

\begin{align}
	\Vc(\chi) := (\chi - \chistar)\T P(\chi - \chistar),\label{eq:Vc}
\end{align}

where $P \in \R^{(\n+N\m) \times (\n+N\m)}$ reads as 

\begin{align}
	P := \begin{bmatrix}
		pI& \frac{1}{N}A\T\1\T  
		\\
		\frac{1}{N}\1 A& pI 
	\end{bmatrix},\label{eq:P}
\end{align}

for some $p > 0$ that we will fix later.

First of all, by the Schur Complement Lemma, we note that for all $p > \bar{p}_0 := \norm{AA\T}$, the matrix $P$ is positive definite.

Moreover, since $\nabla \f$ is $\lipp$-Lipschitz continuous (cf. Assumption~\ref{ass:cost}), once $\reg > 0$ is chosen, it is straightforward to find $\tilde{\lipp} > 0$ such that 

\begin{align}
	\norm{\G(\chi,\chistar)} \leq \tilde{\lipp},
\end{align}

for all $\chi \in \R^{(\n+N\m)}$.

Therefore, the bounds~\eqref{eq:V_quadratic},~\eqref{eq:Lipscitz_Vc}, and~\eqref{eq:nabla_Vc_al} are trivially satisfied.

Hence, we only miss checking the descent condition~\eqref{eq:V_decreasing}.

To this end, we use the definition of $\Vc$ (cf.~\eqref{eq:Vc}) and the reformulation of $\al$ given in~\eqref{eq:G_al} to write

\begin{align}
	&\nabla \Vc(\chi)\T \al(\chi)
	\notag\\
	&= 
	(\chi -\chistar)\T\left(\G(\chi,\chistar)\T P +\right) %P\G(\chi,\chistar)\right)(\chi - \chistar).\label{eq:nabla_W}
\end{align}

Let us focus on the matrix $\G(\chi,\chistar)\T P + P\G(\chi,\chistar)$.

By using the definitions of $\G(\chi,\chistar)$ (cf.~\eqref{eq:G}) and $P$ (cf.~\eqref{eq:P}) and recalling that $\1 \in \ker(\J)$, we have that

\begin{align}
	&\G(\chi,\chistar)\T P + P\G(\chi,\chistar) 
	\notag\\
	&=\begin{bmatrix}
		-\cH(\x,\xstar)& \frac{1}{N}A\T\1\T
		\\
		-\frac{1}{N}\1A& -\reg\J
	\end{bmatrix}\begin{bmatrix}
		pI& A\T\1\T  
		\\
		\1 A& pI 
	\end{bmatrix} 
	\notag\\
	&\hspace{.4cm}
	+ 
	\begin{bmatrix}
		pI& A\T\1\T  
		\\
		\1 A& pI 
	\end{bmatrix}\begin{bmatrix}
		-\cH(\x,\xstar)& -\frac{1}{N}A\T\1\T
		\\
		\frac{1}{N}\1A& -\reg\J
	\end{bmatrix}
	\notag\\
	&=\begin{bmatrix}
		-p\cH(\x,\xstar) + A\T A& -\cH(\x,\xstar)A\T\1\T + \frac{p}{N}A\T\1\T
		\\ 
		-\frac{p}{N}\1A& -\frac{1}{N}\1AA\T\1\T - p\reg\J
	\end{bmatrix}
	\notag\\
	&\hspace{.4cm}
	+ \begin{bmatrix}
		-p\cH(\x,\xstar) + A\T A& -\frac{p}{N}A\T\1\T 
		\\ 
		-\1A\cH(\x,\xstar) + \frac{p}{N}\1A&  -\frac{1}{N}\1AA\T\1\T - p\reg\J
	\end{bmatrix} 
	\notag\\
	&=
	\begin{bmatrix}
		-2p\cH(\x,\xstar) + 2A\T A&  -\cH(\x,\xstar)A\T\1\T
		\\
		-\1A\cH(\x,\xstar)& -\frac{2}{N}\1AA\T\1\T - 2p\reg\J
	\end{bmatrix}.\label{eq:matrix}
\end{align}

Let us focus on the right-down block of the obtained matrix, namely on the term $-\frac{2}{N}\1AA\T\1\T - 2p\reg\J$.

Since $A$ is full-row rank (cf.~Assumption~\ref{ass:full_row}), then $AA\T$ is positive definite.

In turn, since the matrix $\J$ is semi-positive definite and its kernel coincides with the span of $\1$ by construction, the whole term $-\frac{2}{N}\1AA\T\1\T - 2p\reg\J$ is negative definite.

On the other hand, we note that $-2p\cH(\x,\xstar) \leq -2\str I_{n}$ for all $\x\in \R^{\n}$ since $\f$ is strongly convex (cf. Assumption~\ref{ass:cost}).

Therefore, since the parameter $p$ does not appear in the cross terms of the obtained matrix (see~\eqref{eq:matrix}), there exists $\bar{p}_1 > 0$ such that, for all $p > \max\{\bar{p}_0,\bar{p}_1\}$, by Schur Complement Lemma, we have 

\begin{align}\label{eq:quadratic_bound} 
	\G(\chi,\chistar)\T P + P\G(\chi,\chistar) < -  q %I_{(\n+N\m)},
\end{align}

for all $\chi \in \R^{\n}$, where $q := 2\min\{p\str,a,p\reg\}$ in which $a > 0$ denotes the smallest eigenvalue of the matrix $\frac{1}{N}AA\T$.

The result~\eqref{eq:quadratic_bound} allows us to bound~\eqref{eq:nabla_W} as 

\begin{align*}
	\nabla W(\chi)\T \al(\chi) \leq -q\norm{\chi - \chistar}^2,
\end{align*}

for all $\chi \in \R^{(\n+N\m)}$ which concludes the proof.

\bibliography{biblio}

\end{document}